\documentclass[12pt,twoside]{article}

\usepackage{amsgen,amsmath,amstext,amsbsy,amsopn,amsfonts,amssymb}
\usepackage{jltmac2e}

\renewcommand{\version}{} 
\renewcommand{\vol}{??}                             
\renewcommand{\annum}{??}                           
\renewcommand{\title}{\centerline{Leibniz algebras, Lie racks, and digroups}}
\renewcommand{\author}{Michael K. Kinyon}
\renewcommand{\lastname}{Kinyon}

\renewcommand{\entries}{
\nextref{AG}{Andruskiewitsch, N. and M. Gra\~{n}a}
{\textit{From racks to pointed Hopf algebras}}
{Adv. Math. \textbf{178}  (2003) 177--243}

\nextref{Cl}{Clifford, A. H. and G. B. Preston}
{``The Algebraic Theory of Semigroups''}
{Vol. I., Mathematical Surveys \textbf{7},
American Mathematical Society, Providence, R.I. 1961}

\nextref{Cr}{Crans, A.}
{\textit{Lie $2$-algebras}}
{Ph.D. dissertation, University of California Riverside, 2004,
\texttt{www.arxiv.org/abs/math.QA/0409602}}

\nextref{Du}{Duashvili, T.}
{\textit{Witt's theorem for groups with action and free Leibniz algebras}}
{Georgian Math. J. \textbf{11} (2004), 691--712}

\nextref{Fe}{Felipe, R.}
{\textit{Generalized Loday algebras and digroups}}
{preprint, available at
\texttt{www.cimat.mx/reportes/enlinea/I-04-01.pdf}}

\nextref{Fe2}{Felipe, R.}
{\textit{Digroups and their linear representations}}
{East West Math. J., to appear}

\nextref{FR}{Fenn, R. and C. Rourke}
{\textit{Racks and links in codimension two}}
{J. Knot Theory Ramifications \textbf{1} (1992), 343--406}

\nextref{Vienna}{Kinyon, M. K.}
{\textit{The coquecigrue of a Leibniz algebra}}
{presented at \emph{AlanFest},
a conference in honor of the 60th birthday of Alan Weinstein, 
Erwin Schr\"{o}dinger Institute, Vienna, Austria, 4 August 2004.
Available at
\texttt{w3.impa.br/~jair/alanposter/coquecigrue.pdf}}

\nextref{KW}{Kinyon, M. K. and A. Weinstein}
{\textit{Leibniz algebras, Courant algebroids, and multiplications on homogeneous spaces}}
{Amer. J. Math. \textbf{123} (2001) 525--550}

\nextref{KS}{Kosmann-Schwarzbach, Y.}
{\textit{Loday algebra (Leibniz algebra)}}
{Encyclopaedia of Mathematics, Supplement II,
Kluwer Academic Publishers, Dordrecht, 2000, 318--319}

\nextref{KS2}{Kosmann-Schwarzbach, Y.}
{\textit{Derived brackets}}
{Lett. Math. Phys. \textbf{69} (2004), 61--87}

\nextref{Liu}{Liu, K.}
{\textit{A class of grouplike objects}}
{\texttt{www.arXiv.org/math.RA/0311396}}

\nextref{Liu2}{Liu, K.}
{\textit{Transformation digroups}}
{\texttt{www.arXiv.org/math.GR/0409265}}

\nextref{Liu3}{Liu, K.}
{``The Generalizations of Groups''}
{Research Monographs in Math. \textbf{1}, 153 Publishing: Burnaby, 2004.}

\nextref{LoCoq}{Loday, J. L.}
{\textit{Une version non commutative des alg\`{e}bres de Lie:
les alg\`{e}bres de Leibniz}}
{Enseign. Math. \textbf{39} (1993) 269--293}

\nextref{LoDi}{Loday, J. L.}
{\textit{Dialgebras}}
{in ``Dialgebras and related operads'', Lecture Notes in Math. \textbf{1763},
Springer, Berlin, 2001, 7--66}

\nextref{LoPi}{Loday, J. L. and T. Pirashvili}
{\textit{The tensor category of linear maps and Leibniz algebras}}
{Georgian Math. J.  \textbf{5} (1998) 263--276}

\nextref{FO}{Ongay, F.}
{$\varphi$-\textit{dialgebras and a class of matrix ``coquecigrues''}}
{preprint}

\nextref{JD}{Phillips, J. D.}
{\textit{A short basis for the variety of digroups}}
{Semigroup Forum \textbf{70} (2005), 466--470}

\nextref{Sa}{Sabinin, L. V.}
{``Smooth Quasigroups and Loops''}
{Mathematics and its Applications \textbf{492}.
Kluwer Academic Publishers, Dordrecht, 1999}

}

\def\rec{} 
\def\fin{} 

\renewcommand{\address}{Michael K. Kinyon

Department of Mathematics

University of Denver

Denver, CO 80208

\makebox[1in][l]{mkinyon{@}math.du.edu}

\makebox[1in][l]{http://www.math.du.edu/\symbol{126}mkinyon}
}

\newcommand{\seclabel}[1]{\label{sec:#1}} 
\newcommand{\thmlabel}[1]{\label{thm:#1}} 
\newcommand{\lemlabel}[1]{\label{lem:#1}} 
\newcommand{\corolabel}[1]{\label{coro:#1}} 
\newcommand{\deflabel}[1]{\label{def:#1}} 
\newcommand{\exlabel}[1]{\label{ex:#1}} 
\newcommand{\eqlabel}[1]{\label{eq:#1}} 

\newcommand{\thmref}[1]{\ref{thm:#1}} 
\newcommand{\lemref}[1]{\ref{lem:#1}} 
\newcommand{\exref}[1]{\ref{ex:#1}} 
\renewcommand{\eqref}[1]{\ref{eq:#1}} 
\newcommand{\peqref}[1]{(\eqref{#1})} 

\newcommand{\iv}{^{-1}}
\newcommand{\lp}{\vdash}
\newcommand{\rp}{\dashv}
\newcommand{\lie}[2]{\lbrack #1 , #2 \rbrack}
\newcommand{\fg}{\mathfrak{g}}
\newcommand{\fh}{\mathfrak{h}}
\newcommand{\cE}{\mathcal{E}}
\newcommand{\cA}{\mathcal{A}}
\newcommand{\cS}{\mathcal{S}}

\newcommand{\End}{\mathrm{End}}
\newcommand{\Aut}{\mathrm{Aut}}
\newcommand{\Der}{\mathrm{Der}}
\newcommand{\Ad}{\mathrm{Ad}}
\newcommand{\ad}{\mathrm{ad}}

\begin{document}
\firstpage

\begin{abstract}
The ``coquecigrue" problem for Leibniz algebras is that of finding
an appropriate generalization of Lie's third theorem, that is, of
finding a generalization of the notion of group such that Leibniz
algebras are the corresponding tangent algebra structures.
The difficulty is determining exactly what properties this generalization
should have. Here we show that \emph{Lie racks}, smooth
left distributive structures, have Leibniz algebra structures on their
tangent spaces at certain distinguished points. One way of producing
racks is by conjugation in \emph{digroups}, a generalization of group
which is essentially due to Loday. Using semigroup theory, we show that
every digroup is a product of a group and a trivial digroup. We partially solve the 
coquecigrue problem by showing that to each Leibniz algebra that
splits over an ideal containing its ideal generated by squares,
there exists a special type of Lie digroup with tangent algebra
isomorphic to the given Leibniz algebra. The general coquecigrue
problem remains open, but Lie racks seem to be a promising direction.
\end{abstract}

\section{Introduction}
\seclabel{intro}

A \emph{Leibniz algebra} $(\fg, \lie{\cdot}{\cdot})$ is a vector
space $\fg$ together with a bilinear mapping
$\lie{\cdot}{\cdot} : \fg\times \fg\to \fg$ satisfying the derivation
form of the Jacobi identity
\begin{equation}
\eqlabel{Jacobi}
\lie{X}{\lie{Y}{Z}} = \lie{\lie{X}{Y}}{Z} + \lie{Y}{\lie{X}{Z}}
\end{equation}
for all $X,Y,Z\in \fg$. Every Lie algebra is a Leibniz algebra, but
the bracket in a Leibniz algebra need not be skew-symmetric.

One of the outstanding problems in the theory of Leibniz algebras
is that of finding an appropriate generalization of Lie's ``third theorem'',
which associates a (local) Lie group to any (real or complex) Lie algebra.
The most challenging aspect of the problem is to determine what should be the
correct generalization of the notion of group. So little is known about
what properties these group-like objects should have that Loday dubbed
them ``coquecigrues" \cite{LoCoq}. The author and A. Weinstein made
an attempt to understand coquecigrues in \cite{KW} by using smooth,
nonassociative multiplications defined on reductive homogeneous spaces
that are associated to Leibniz algebras. However, as already discussed in
that paper, the objects studied therein are not a satisfactory solution to the
coquecigrue problem, because in the case where the given Leibniz algebra is a Lie
algebra, the object associated to it in \cite{KW} is not the Lie group. At the
very least, the correct generalization of Lie's third theorem to Leibniz 
algebras should reduce to the usual theorem for Lie algebras; that is, the
coquecigrue of a Lie algebra should be a Lie group.

This paper offers a different partial solution to the coquecigrue problem.
By dissecting one of the proofs that the tangent space at the unit
element of a Lie group is a Lie algebra, we find that the essential
properties leading to the Jacobi identity are not so much encoded in
the group multiplication as they are in conjugation. This is hardly an
original observation. For instance, recently it has been made in a
categorical context by Crans, who is interested in the similar problem
of the relationship between Lie $2$-groups and Lie $2$-algebras \cite{Cr}. In
any case, the observation that conjugation is what induces the Jacobi
identity leads us to
consider \emph{Lie racks}, that is, (pointed) manifolds with a smooth, left
distributive binary operation. We propose that whatever the correct notion of
coquecigrue might turn out to be, the Leibniz algebra of a coquecigrue
should be obtained from differentiating a Lie rack structure associated
to the coquecigrue. Our primary example (Example \exref{lie-rack}) of
a nongroup Lie rack is a product of a Lie group and a vector space upon
which the group acts; we call these linear Lie racks.

The notion of Lie rack itself might seem to be a reasonable candidate for
coquecigrue, but this is premature, because it turns out that every
Leibniz algebra can be realized as the tangent Leibniz algebra of some rack,
but in the Lie algebra case, the rack in question is not the conjugation
rack of the Lie group. However, there are other rack structures that are
more promising, at least in the split case. In Theorem \thmref{split-rack},
we show that the Leibniz algebra
$\fg$ associated to a linear Lie rack splits over its ideal $\cE$ generated
by squares, that is, there exists a Lie algebra $\fh \subset \fg$ such that
$\fg = \cE \oplus \fh$, a direct sum of vector spaces. In this case,
$\fg$ is the demisemidirect product of $\cE$ with $\fh$ in the sense
of \cite{KW}. Conversely, if $\fg$ is a split Leibniz algebra, then
there exists a linear Lie rack with $\fg$ as its tangent Leibniz algebra.

There are algebraic structures that are intermediate between racks
and groups considered as racks. One such structure, which is suggested
by the work of Loday \cite{LoDi}, is the notion of \emph{digroup}.
We will give an axiomatic description of digroups in \S{4},
but essentially, they are sets $G$ with two binary operations $\lp$ and $\rp$
so that, in the jargon of semigroup theory \cite{Cl},  $(G,\lp)$ is a right group,
$(G,\rp)$ is a left group, and there are some compatibility conditions
ensuring that the sets of unit elements of the two structures
coincide. Using basic semigroup theory, we show that every digroup
is a product of a group and a ``trivial" digroup.

Conjugation in groups generalizes to digroups so that every digroup
has an associated rack structure. In the smooth case, then, every
Lie digroup $G$ (that is, digroup in which $\lp$ and $\rp$ are smooth
operations with respect to an underlying manifold structure) has an associated
Leibniz algebra $\fg = T_1G$, where $1$ is a distinguished unit element.
The rack structure associated to a linear digroup is a linear Lie rack,
and conversely, every linear Lie rack is induced by a linear digroup
structure. As a corollary, we obtain Lie's third theorem for
split Leibniz algebras: to each split Leibniz algebra $\fg$, there exists a
Lie digroup $G$ such that $T_1G$ is a Leibniz algebra isomorphic to $\fg$.

Within a span of a year, the notion of digroup discussed in this paper
has been introduced independently three times. The first appearance
of the definition and a discussion of basic properties seems to have been
that of the author \cite{Vienna} in a conference talk on 4 August 2004.
Next appeared a preprint of K. Liu \cite{Liu}, and then
another by R. Felipe \cite{Fe}. In each case, not only were the definitions
of digroup nearly identical, but so was the choice of the term ``digroup" to denote
the concept. This is not as surprising as it might seem: the definition of
dimonoid (which probably should have been called disemigroup) and the
``di-" terminology are already contained in the work of J. Loday \cite{LoDi}.
Digroups are just dimonoids in which every element is invertible in an appropriate
sense. Thus one should really see the notion of digroup as being already
implicit in Loday's work. The observation that a digroup is a right group and
a left group with a common set of unit elements seems to be new to this paper.
Additional recent papers on digroups and/or coquecigrues are
\cite{Fe2, Liu2, FO}, and Liu has published a monograph \cite{Liu3},
which also contains a wider notion of digroup than that studied here.

Worth mentioning is the interesting work of T. Datuashvili \cite{Du}, who has
approached the coquecigrue problem from the point of view of generalizing the
Magnus-Witt functor from the category of groups to the category of Lie algebras.

I would like to thank Alan Weinstein for introducing me to the problem of
``integrating" Leibniz algebras. After witnessing \cite{Vienna},
he made a couple of key remarks which greatly helped to simplify the presentation
of the notion of digroup. He also made very useful comments on an early version
of this paper. I would also like to thank Yvette Kosmann-Schwarzbach 
and Dmitri Roytenberg for their encouragement.

\section{Leibniz algebras}
\seclabel{leibniz}

This brief section summarizes basic facts about Leibniz algebras, and defines
split Leibniz algebras. In this paper, the underlying field of any vector space
will be the real or complex numbers. The problem of ``integrating" Leibniz algebras
over other fields to an appropriate generalization of the notion of algebraic group
is an interesting one, but will not be addressed here.

Let $\fg$ be a Leibniz algebra as defined in \S{1}.
If we define the left multiplication map by $\ad(X)Y := \lie{X}{Y}$,
then the Jacobi identity \peqref{Jacobi} can be summarized by the
assertion that $\ad(X) \in \Der(\fg)$ for all $X\in \fg$, where
$\Der(\fg)$ denotes the Lie algebra of derivations of $\fg$.

It should be noted the convention adopted here for the Jacobi
identity \peqref{Jacobi} defines what is known as a \emph{left}
Leibniz algebra. The opposite algebra of a left Leibniz algebra
is a right Leibniz algebra, in which the right multiplication
maps are derivations. It is somewhat more common in the
literature to work with right Leibniz algebras; however, left Leibniz algebras seem a bit closer
to the spirit of the coquecigrue problem, as will be seen in \S{3}.
Leibniz algebras are also known as Loday algebras \cite{KS},
because they were introduced and popularized by Loday \cite{LoCoq}.

Denote by $\cS := \langle \lie{x}{x} : x\in \fg \rangle$ the ideal of $\fg$
generated by all squares. Then $\cS$ is the minimal ideal with respect 
to the property that $\fh := \fg / \cS$ is a Lie algebra. The quotient mapping
$\pi : \fg \to \fh$ is a homomorphism of Leibniz algebras, that is,
$\pi(\lie{X}{Y}) = \lie{\pi(X)}{\pi(Y)}$ for all $X,Y\in \fg$.
It is also a morphism of left $\fh$-modules, that is,
$\pi (\xi X) = \lie{\xi}{\pi(X)}$ for all $X\in \fg$, $\xi\in \fh$.

Conversely, Let $\fh$ be a Lie algebra, $\fg$ an $\fh$-module, and
$\pi : \fg \to \fh$ an $\fh$-module morphism. Assume without loss that
$\pi$ is an epimorphism, and define the structure of a Leibniz algebra on
$\fg$ by $\lie{X}{Y} := \pi(X) Y$. Then $\pi$ becomes an epimorphism
of Leibniz  algebras. To view Leibniz algebras in this way is to regard them
as being Lie algebra objects in the infinitesimal tensor category of linear maps;
see \cite{LoPi} for further details. Incidentally, this particular view of
Leibniz algebras may be the key to relating the study of Leibniz algebras
to Crans' study of Lie $2$-algebras, since the latter can be viewed in a
similar fashion \cite{Cr}.

Let $\fg$ be a Leibniz algebra, and let $\cE \subset \fg$ be an ideal such that
$\cS \subseteq \cE \subseteq \ker(\ad)$. 
Then $\fg$ \emph{splits} over $\cE$
if there exists a Lie subalgebra $\fh \subset \fg$ such that
$\fg = \cE \oplus \fh$, a direct sum of vector spaces. In this
case, for $u,v\in \cE$, $X,Y\in \fh$, we have
\begin{equation}
\eqlabel{split-bracket}
\lie{u+X}{v+Y} = Xv + \lie{X}{Y},
\end{equation}
since $\cE \in \ker(\ad)$, that is, $\lie{u}{\cdot} \equiv 0$
for all $u\in \cE$.

Conversely, given a Lie algebra $\fh$ and an $\fh$-module $V$,
set $\fg := V\oplus \fg$, and define
a bracket on $\fg$ by \peqref{split-bracket} for $u,v\in V$,
$X,Y\in \fh$. Then $(\fg, \lie{\cdot}{\cdot})$ is a Leibniz algebra
called the \emph{demisemidirect product} of $V$ and $\fh$ \cite{KW}.
The ideal of $\fg$ generated by squares is $\cS \cong \fh V$.
The kernel of $\ad$ is
$\ker(\ad) = V\oplus \{ X \in \fh : Xv = 0 \ \forall v\in V\} 
\cap Z(\fh)$, where $Z(\fh)$ denotes the center of $\fh$.
Finally, $V \cong V\oplus \{0\}$ is an ideal such that
$\cS \subseteq V \subseteq \ker(\ad)$. We have $\fg / V \cong \fh$,
and so $\fg$ splits over $V$.

We summarize this discussion as follows.

\begin{theorem}
\thmlabel{split-Leibniz}
Let $\fg$ be a Leibniz algebra, and let $\cE$ be an ideal of
$\fg$ such that $\cS \subseteq \cE \subseteq \ker(\ad)$.
Then $\fg$ splits over $\cE$ if and only if $\fg$ is a demisemidirect
product of $\cE$ with a Lie algebra $\fh$.
\end{theorem}

Note that if a Leibniz algebra $\fg$ splits over $\cS$
itself with complementary Lie subalgebra $\fh$, then by
\peqref{split-bracket}, $\fh \cdot \cS = \cS$. Conversely,
if $\fg = V\oplus \fh$ is a demisemidirect product and
if $\fh \cdot V = V$, then $\cS \cong V$.

A Leibniz algebra may split over more than one ideal.

\begin{example}
\exlabel{split1}
Let $V := \mathbb{R}^n$ and on
$\fg := V \oplus gl(V)\oplus gl(V)$, define
$$
\lie{u + X + Y}{v + U + V} := Yv + \lie{X}{U} + \lie{Y}{V} 
$$
for $u,v\in V$, $X,Y,U,V\in gl(V)$. Then $\cS \cong V$ and
$\ker(\ad) = V \oplus \{ aI : a\in \mathbb{R}\} \oplus \{ 0 \}$,
where $I\in gl(V)$ denotes the identity matrix. Here $\fg$
splits over $\cS$ with complement $gl(V)\oplus gl(V)$, and
$\fg$ also splits over $\ker(\ad)$ with complement
$sl(V)\oplus gl(V)$.
\end{example}

If $\fg$ is a Lie algebra which does not split over its center
$Z(\fg)$, then $\fg$ is obviously
a Leibniz algebra that splits (trivially) over its
ideal generated by squares, but not over $\ker(\ad)
= \ker(\mathrm{ad}) = Z(\fg)$. On the other hand, 
a Leibniz algebra may split over $\ker(\ad)$ without
splitting over its ideal generated by squares.

\begin{example}
\exlabel{split2}
Let $V := \mathbb{R}^2$, let
$\fh := \displaystyle\{
\begin{pmatrix} 0 & a \\ 0 & 0 \end{pmatrix} : a\in \mathbb{R}\}$,
and let $\fg = V\oplus \fh$ be the demisemidirect product. 
Then $\ker(\ad) \cong V$ and 
$\cS \cong \mathbb{R}\cdot \begin{pmatrix} 1 \\ 0 \end{pmatrix}$.
While $\fg$ splits over $\ker(\ad)$, it is easy to check that
$\fg$ does not split over $\cS$.
\end{example}

To conclude this section, we will briefly discuss the notion of
dialgebra \cite{LoDi}, as it motivates the notion of digroup.
The commutator bracket in an associative algebra gives the underlying vector
space the structure of a Lie algebra. One way to generalize
this idea so as to obtain Leibniz algebra brackets that are
not skew-symmetric is to use two different associative algebra
products.

\begin{Definition}
\deflabel{dialgebra}
A \emph{dialgebra} $(\cA,\lp,\rp)$ is a vector space $\cA$
together with two associative, bilinear mappings 
$\lp, \rp : \cA\times \cA\to \cA$ satisfying the following axioms.
For all $x,y,z\in \cA$,
\begin{enumerate}
\item[(D1)] $x\lp (y\rp z) = (x\lp y)\rp z$
\item[(D2)] $x\rp (y\lp z) = x\rp (y\rp z)$
\item[(D3)] $(x\rp y)\lp z = (x\lp y)\lp z$
\end{enumerate}
\end{Definition}

The axioms (D1)-(D3) will appear again in \S{4},
and so we postpone a detailed discussion until then.

Given a dialgebra $(\cA,\lp,\rp)$, defining a bracket
by $\lie{x}{y} := x \lp y - y \rp x$ turns
$(\cA,\lie{\cdot}{\cdot})$ into a Leibniz algebra.
Note that our use of $\lp$ and $\rp$ in this bracket
is the opposite of that of Loday. This convention matches our
preference for left Leibniz algebras instead of right Leibniz algebras.

\begin{example}
\exlabel{dialgebra}
Let $V$ be a vector space. On $\cA := V \oplus \End(V)$, define
\begin{align*}
(u,X) \lp (v,Y) &:= (Xv, XY) \\
(u,X) \rp (v,Y) &:= (0,XY)
\end{align*}
Then the associated Leibniz algebra $(\cA, \lie{\cdot}{\cdot})$ is exactly the
demisemidirect product of $V$ with the Lie algebra $gl(V) = \End(V)$.
\end{example}

\section{Lie racks}
\seclabel{racks}

Pointed racks are algebraic structures that encapsulate some
of the properties of group conjugation.

\begin{Definition}
\deflabel{rack}
A \emph{pointed rack} $(Q,\circ,1)$ is a set $Q$ with a binary
operation $\circ$ and a distinguished element $1\in Q$ such that
the following axioms are satisfied.
\begin{enumerate}
\item $x\circ (y\circ z) = (x\circ y) \circ (x\circ z)$ (left distributivity)
\item For each $a,b\in Q$, there exists a unique $x\in Q$
such that $a\circ x = b$
\item $1\circ x = x$ and $x\circ 1 = 1$ for all $x\in Q$.
\end{enumerate}
\end{Definition}

For a magma $(Q,\circ)$, let $\Aut(Q)$ denote the set of permutations
of $Q$ preserving $\circ$, that is, a permutation $\psi$ of $Q$ is in $\Aut(Q)$ if
and only if $\psi (x\circ y) = \psi(x)\circ \psi(y)$ for all $x,y\in Q$.
If we denote the left translations in a rack $(Q,\circ)$ by
$\phi(x)y := x\circ y$, then the left distributive axiom simply asserts
that $\phi(x) \in \Aut(Q)$ for all $x\in Q$.

In a group $G$, the operation $x\circ y := xyx\iv$ makes $(G,\circ,1)$
into a pointed rack, where $1$ is the identity element of $G$.
Left (or right) distributive structures have been studied under
an overabundance of names in the literature. ``Rack" seems to
be the most fashionable name in the case where each $\phi(x)$
is bijective. A good survey of racks and more specialized structures,
such as quandles and crossed sets, can be found in \cite{AG}.
In general, racks need not be pointed, but all racks considered
in this paper are. 

We are primarily interested in racks in which the algebraic
structure is compatible with an underlying manifold structure.
In this paper, ``smooth manifold" means $C^{\infty}$-smooth.

\begin{Definition}
\deflabel{lierack}
A \emph{Lie rack} $(Q,\circ,1)$ is a smooth manifold $Q$
with the structure of a pointed rack such that the rack
operation $\circ: Q\times Q\to Q$ is a smooth mapping.
\end{Definition}

The conjugation operation in a Lie group makes it into a Lie rack.
Here is another example.

\begin{example}
\exlabel{lie-rack}
Let $H$ be a Lie group and let $V$ be an $H$-module.
On $Q := V\times H$, define a binary operation $\circ$ by
\begin{equation}
\eqlabel{lie-rack}
(u, A)\circ (v, B) := (Av, ABA\iv )
\end{equation}
for all $u,v\in V$, $A,B\in H$. Setting $\mathbf{1} = (0,1)$,
we have that $(Q,\circ,\mathbf{1})$ is a Lie rack, which we
call a \emph{linear Lie rack}.
\end{example}

By defining a Lie rack to be pointed, we do not wish to suggest
that smooth unpointed racks are of no interest. Smooth left
distributive structures have certainly been studied in the literature;
see the bibliography of \cite{Sa}. However, our Lie racks $Q$ have a distinguished tangent
space $T_1Q$ with an algebraic structure of its own. Our approach to
that structure is modeled upon one of the routes to the Lie algebra of a Lie group.
In this approach, the idea is to differentiate the conjugation
operation to obtain the adjoint representation of the group,
and then differentiate again to obtain a mapping which is used
to define the Lie bracket, and which then becomes the adjoint
representation of the Lie algebra. 

Suppose now that $Q$ is a Lie rack. For each $x\in Q$,
$\phi(x)1 = 1$, and so we may apply the tangent functor $T_1$
to $\phi(x) : Q\to Q$ to obtain a linear mapping
$\Phi(x) := T_1\phi(x) : T_1Q \to T_1Q$. Since each $\phi(x)$ is invertible, 
we have each $\Phi(x)\in GL(T_1Q)$. Now the
mapping $\Phi : Q \to GL(T_1Q)$ satisfies $\Phi(1) = I$, where
$I\in GL(T_1Q)$ is the identity mapping. Thus we 
may differentiate again to obtain a mapping
$\ad : T_1Q \to gl(T_1Q)$. 
Here we are making the usual identification of the tangent
space at the identity element of $GL(V)$ for a vector space
$V$ with the general linear Lie algebra $gl(V)$. Now we set
\begin{equation}
\eqlabel{bracket}
\lie{X}{Y} := \ad(X)Y
\end{equation}
for all $X,Y\in T_1Q$.

In terms of the left multiplications $\phi(x)$, the left distributive
property of racks can be expressed by the equation
\begin{equation}
\eqlabel{distrib}
\phi(x)\phi(y) z = \phi(\phi(x)y)\phi(x)z .
\end{equation}
We differentiate \peqref{distrib} at $1\in Q$, first with respect to $z$ and then with respect
to $y$ to obtain
\begin{equation}
\eqlabel{Phi-aut}
\Phi(x) \lie{Y}{Z} = \lie{\Phi(x) Y}{\Phi(x) Z}
\end{equation}
for all $x\in Q$, $Y,Z\in T_1Q$.
This expresses the condition that for each $x\in Q$,
$\Phi(x) \in \Aut(T_1Q, \lie{\cdot}{\cdot})$.

Next we differentiate \peqref{Phi-aut} at $1$ to get
\begin{equation}
\eqlabel{tangent-Leibniz}
\lie{X}{\lie{Y}{Z}} = \lie{\lie{X}{Y}}{Z} + \lie{Y}{\lie{X}{Z}}
\end{equation}
for all $X,Y,Z\in T_1Q$.

Summarizing, we have shown the following.

\begin{theorem}
\thmlabel{rack-Leibniz}
Let $(Q,\circ,1)$ be a Lie rack, and let $\fg = T_1Q$.
Then there exists a bilinear mapping
$\lie{\cdot}{\cdot} : \fg\times \fg\to \fg$
such that
\begin{enumerate}
\item $(\fg, \lie{\cdot}{\cdot})$ is a Leibniz algebra,
\item For each $x\in Q$, the tangent mapping $\Phi(x) = T_1\phi(x)$
is an automorphism of $(\fg,\lie{\cdot}{\cdot})$, 
\item if $\ad : \fg \to gl(\fg)$ is defined by
$X\mapsto \ad(X) = (Y\mapsto \lie{X}{Y} )$, then $\ad = T_1 \Phi$.
\end{enumerate}
\end{theorem}

We will refer to the Leibniz algebra structure on the tangent space
at the distinguished element $1$ of a Lie rack as being the tangent
Leibniz algebra of the rack. To illustrate Theorem \thmref{rack-Leibniz},
let us explicitly describe the tangent Leibniz algebra structure for the linear Lie rack
$(Q,\circ,\mathbf{1})$ of Example \exref{lie-rack}. Let $\fh$ be the
Lie algebra of the Lie group $H$. Then we may identify 
$T_1Q$ with $\fg := V\oplus \fh$. For $u\in V$, $A\in H$, the 
tangent mapping $\Phi(u,A) = T_1 \phi(u,A) : \fg \to \fg$ is given by
\begin{equation}
\eqlabel{tanmap}
\Phi(u,A) (v + X) = Av + \Ad(A) X
\end{equation}
for all $v\in V$, $X\in \fh$, where $\Ad : H \to GL(\fh)$ is the
adjoint representation. Differentiating this, we find that
\begin{equation}
\eqlabel{tanLeibniz}
\lie{u + X}{v + Y} = Xv + \lie{X}{Y}
\end{equation}
for all $u,v\in V$, $X,Y\in \fh$. Comparing with \peqref{split-bracket}, we 
see that the tangent Leibniz algebra for the Lie rack $Q = V\times H$ of
Example \exref{lie-rack} is exactly the demisemidirect product of $V$ with $\fh$.

So far, we have shown one direction of the following.

\begin{theorem}
\thmlabel{split-rack}
Let $H$ be a Lie group with Lie algebra $\fh$, let $V$ be an $H$-module,
and let $(Q, \circ, \mathbf{1})$ be the linear Lie rack defined by \peqref{lie-rack},
where $Q = V\times H$. Then the tangent Leibniz algebra of $Q$ is the demisemidirect
product $\fg = V \oplus \fh$ with bracket given by \peqref{tanLeibniz}.

Conversely, let $\fg$ be a split Leibniz algebra. Then there exists a
linear Lie rack $Q$ with tangent Leibniz algebra isomorphic to $\fg$.
\end{theorem}

\begin{proof}
Only the second assertion remains to be shown. Let $\fg = \cE \oplus \fh$ be a splitting
of $\fg$, where $\cE$ is an ideal with $\cS \subseteq \cE \subseteq \ker(\ad)$
and $\fh$ is a Lie subalgebra.
Recall that $\fg$ is then a demisemidirect product of $\cE$ with $\fh$. Let
$H$ be a connected Lie group with Lie algebra $\fh$. Set $Q = \cE\times H$, and
note that we may identity $\fg$ with $T_1Q$. Give $Q$ the Lie rack structure
$(Q,\circ, \mathbf{1})$ where $\circ$ is given by \peqref{lie-rack}. Then the
result follows from the discussion preceding the statement of the theorem.
\end{proof}

It is natural to speculate that Lie racks themselves provide an answer to the
coquecigrue problem. That this is not the case is evidenced by the observation
that \emph{every} Leibniz algebra is the tangent Leibniz algebra of a rack.
Indeed, let $\fg$ be a Leibniz algebra with $\ad(X)Y := \lie{X}{Y}$.
On $\fg$, define
$$
X\circ Y := \exp \ad(X) Y
$$
Then
$$
\begin{array}{rcl}
X\circ (Y\circ Z) &=& \exp{\ad(X)}\exp{\ad(Y)}Z \\
&=& \exp (\exp \ad(X)\ad(Y)) \exp{\ad(X)}Z \\
&=& \exp(\ad(\exp{\ad(X)Y}))\exp{\ad(X)}Z \\
&=& (X\circ Y)\circ (X\circ Z)
\end{array}
$$
Also $X\circ 0 = 0$ and $0\circ X = X$. Thus $(\fg, \circ)$ is
a Lie rack. It is easy to check that the tangent Leibniz algebra
is $\fg$ itself.

In case $\fg$ is a Lie algebra, the rack structure $(\fg,\circ)$ was
first noted by Fenn and Rourke \cite{FR}. Since the coquecigrue of
a Lie algebra is supposed to be a Lie group, this rack does not meet
a basic requirement for being the coquecigrue of a Leibniz algebra.

In the Lie algebra case, it is natural to ask how
$(\fg,\circ)$ is related to the associated conjugation rack of a Lie
group $G$ for $\fg$. The relationship can be seen in the Lie group structure
on the tangent bundle $TG = \fg \rtimes G$. Conjugation in
$TG$ is given by
$$
(X,a)\circ (Y,b) = (X + \Ad(a)Y - \Ad(a\iv)X, aba\iv)
$$
for $X,Y\in \fg$, $a,b\in G$. Now consider the graph
$\{ (X, \exp X) : X \in \fg \}$ of the exponential mapping.
This graph is not, in general, a subgroup of $TG$. However, it is
a subrack; from the above we have
$$
\begin{array}{rcl}
(X, \exp X)\circ (Y, \exp Y) &=& ( \Ad(\exp X)Y , \exp X \exp Y \exp (-X)) \\
&=& (\Ad(\exp X)Y, \exp(\Ad(\exp X)Y) )
\end{array}
$$
This subrack is obviously just a copy of $(\fg,\circ)$. In fact, it is
the graph of the rack homomorphism
$(\fg,\circ) \to (\Aut(\fg), \circ) ; X\mapsto \Ad(\exp X)$.

Since Lie racks are not a sufficient answer to the coquecigrue problem,
more structure must be needed which induces an associated rack structure.
In the next section, we study a candidate for this additional structure.

\section{Digroups}
\seclabel{digroups}

Just as dialgebras with two distinct associative operations lead to nonLie
Leibniz algebras through a generalized notion of bracket, so will we obtain
Lie racks that are not groups via a generalized notion of conjugation.

\begin{Definition}
\deflabel{digroup}
A \emph{disemigroup} $(G,\lp,\rp)$ is a set $G$ together
with two binary operations $\lp$ and $\rp$ satisfying the
following axioms. For all $x,y,z\in G$,
\begin{enumerate}
\item[(G1)] $(G,\lp)$ and $(G,\rp)$ are semigroups
\item[(G2)] $x\lp (y\rp z) = (x\lp y)\rp z$
\item[(G3)] $x\rp (y\lp z) = x\rp (y\rp z)$
\item[(G4)] $(x\rp y)\lp z = (x\lp y)\lp z$
\end{enumerate}
A disemigroup is a \emph{dimonoid} if
\begin{enumerate}
\item[(G5)] $\exists 1\in G$ such that $1\lp x = x\rp 1 = x$ for all $x\in G$.
\end{enumerate}
A dimonoid is a \emph{digroup} if
\begin{enumerate}
\item[(G6)] $\forall x\in G, \exists x\iv \in G$ such that $x\lp x\iv = x\iv \rp x = 1$.
\end{enumerate}
\end{Definition}

Loday used the term ``dimonoid" to refer to what we have
called a disemigroup \cite{LoDi}. We have made a slight change in the
terminology to be more consistent with standard usage in
semigroup theory. 
An element $e$ in a disemigroup is called a \emph{bar-unit}
if it satisfies $e\lp x = x\rp e = x$ for all $x$. Axiom (G5)
asserts that a bar-unit exists in a dimonoid, but it is not assumed to be
unique. Note that a digroup is a group if and only if $\lp\  =\  \rp$
if and only if $1$ is the unique bar-unit.
As noted in \S{1}, the notion of
digroup has recently been introduced independently three times
\cite{Fe,Vienna,Liu}, but one should really think of it as being already
implicit in Loday's work.

We will use the associativity of the operations $\lp$ and $\rp$
as well as axiom (G2) to drop parentheses from expressions
whenever possible. Thus, for instance, $x\lp y\rp z$ is 
unambiguous by (G2).

J.D. Phillips has recently shown that in the presence of axioms
(G5) and (G6), axioms (G3) and (G4) are redundant, and in fact, 
axioms (G1), (G2), (G5), (G6) are an independent set of axioms
for digroups \cite{JD}.

\begin{example}
\exlabel{easy}
Let $H$ be a group and $M$ a set on which $H$ acts on the
left. Suppose there exists a fixed point $e\in M$, that is, $he = e$
for all $h\in H$, and suppose that $H$ acts transitively on 
$M\backslash \{e\}$. On $G := M\times H$, define
\begin{eqnarray}
(u,h) \lp (v,k) &= (hv, hk) \eqlabel{lp1}\\
(u,h) \rp (v,k) &= (u, hk) \eqlabel{rp2}
\end{eqnarray}
for all $u,v\in M$, $h,k\in H$. Then $(G,\lp,\rp)$ is a digroup
with distinguished bar-unit $(e,1)$. The inverse of $(u,h)$
is $(e,h\iv)$.
\end{example}

We will see later that, in a sense to be made more precise,
every digroup is of the type described in 
Example \exref{easy}.

Understanding the definition of digroup becomes easier if we draw
upon a bit of semigroup theory. Let $(G,\lp,\rp)$ be a digroup.
The semigroup $(G,\lp)$ is presumed to have a left neutral
element $1$ and right inverses $x\iv$ for each $x\in G$.
Semigroups with this additional structure are called 
\emph{right groups} \cite{Cl}. (This is not how right groups
are usually defined in the semigroup literature, but is instead
a characterization.) Similarly, $(G,\rp)$, which
has a right neutral element $1$ and left inverses $x\iv$ for
each $x\in G$, is a \emph{left group}. 

We collect some basic facts about right groups in the next
result; see, e.g., \cite{Cl}.

\begin{lemma}
\lemlabel{left-groups}
Let $(G,\lp)$ be a right group, let $J := \{ x\iv : x\in G \}$, and
let $E := \{ e : e\lp x = x \}$
\begin{enumerate}
\item $x\iv \lp x \lp y = x \lp x\iv \lp y = y$ for all $x,y\in G$.
\item $x\lp 1 = (x\iv)\iv$ for all $x\in G$.
\item $((x\iv)\iv)\iv = x\iv$ for all $x\in G$.
\item $(x\lp y)\iv = y\iv \lp x\iv$ for all $x,y\in G$.
\item $J$ is a group and $E$ is a right zero semigroup.
\item $G\to J; x\mapsto (x\iv)\iv$ is an epimorphism of right groups
with kernel $E$.
\item $G\to E; x\mapsto x\iv \lp x$ is an epimorphism of right groups
with kernel $J$.
\item $G = J\lp E$ is isomorphic to the direct sum of $J$
and $E$.
\end{enumerate}
\end{lemma}

The set $E$ is sometimes called the ``halo'' in papers on dialgebras
\cite{LoDi}.

When we need results about left groups corresponding to
those in Lemma \lemref{left-groups} in the sequel, we will
simply refer to the ``left group dual" of the appropriate assertion.

Now let $(G,\lp,\rp)$ be a digroup. In order to interpret
the digroup axioms, we introduce mappings on $G$ as follows:
$$
L_{\lp}(x) y = R_{\lp}(y) x := x \lp y \qquad \text{and} \qquad
L_{\rp}(x) y = R_{\rp}(y) x := x \rp y
$$
Since $(G,\lp)$ and $(G,\rp)$ are semigroups, we have the usual  relations
$$
\begin{array}{rclcrcl}
L_{\lp}(x) R_{\lp}(y) &=& R_{\lp}(y) L_{\lp}(x)
& & 
L_{\rp}(x) R_{\rp}(y) &=& R_{\rp}(y) L_{\rp}(x)
\\
L_{\lp}(x\lp y) &=& L_{\lp}(x) L_{\lp}(y)
& \qquad \text{and} \qquad &
L_{\rp}(x\rp y) &=& L_{\rp}(x) L_{\rp}(y)
\\
R_{\lp}(x\lp y) &=& R_{\lp}(y) R_{\lp}(x)
& &
R_{\rp}(x\rp y) &=& R_{\rp}(y) R_{\rp}(x)
\end{array}
$$
If we let $G^G$ denote the semigroup of all mappings on $G$, then
$x\mapsto L_{\lp}(x)$ and $x\mapsto R_{\lp}(x)$ are, respectively,
a homomorphism and an antihomomorphism from $(G,\lp)$ to $G^G$.
Also, $x\mapsto L_{\rp}(x)$ and $x\mapsto R_{\rp}(x)$ are, respectively,
a homomorphism and an antihomomorphism from $(G,\rp)$ to $G^G$.
In particular, $L_{\lp}(1) = R_{\lp}(1) = L_{\rp}(1) = R_{\rp}(1) = I$.
Since $(G,\lp)$ is a right group, the image of $L_{\lp}$ lies in $G!$, the
symmetric group on $G$. In particular, by Lemma \lemref{left-groups}(1),
$L_{\lp}(x\iv) = L_{\lp}(x)\iv$. Similarly, since $(G,\rp)$ is a left
group, the image of $R_{\rp}$ lies in $G!$, and
$R_{\rp}(x\iv) = R_{\rp}(x)\iv$.

Turning to the compatibility axioms, (G2) reads
$$
L_{\lp}(x) R_{\rp}(z) = R_{\lp}(z) L_{\rp}(x)
$$
Axioms (G3) and (G4) read, respectively
$$
R_{\rp}(y \rp z) = R_{\rp}(y \lp z)
\qquad \text{and} \qquad
L_{\lp}(x \lp y) = L_{\lp}(x \lp y)
$$
So the mapping $L_{\lp} : G\to G!$, in addition to be a homomorphism
from $(G,\lp)$ to $G!$, is also a homomorphism from $(G,\rp)$ to $G!$.
Similarly, $R_{\rp} : G\to G!$, in addition to being an antihomomorphism
from $(G,\rp)$ to $G!$ is also an antihomomorphism from $(G,\lp)$ to $G!$.

Liu \cite{Liu} has shown that the subset $L_{\lp}(G) \times L_{\rp}(G)
\subset G!\times G!$ has a natural digroup structure containing as a
subdigroup a copy of $(G,\lp,\rp)$. He interprets this to be a type of Cayley
representation of digroups. Here we will obtain a different representation
using basic semigroup theory.

The kernel of the homomorphism $L_{\lp} : G\to G!$ is the set of
all bar-units of $G$, that is, the set
\begin{equation}
\eqlabel{bar-units}
E := \{ e \in G : e\lp x = x\rp e = x \ \ \forall x\in G \} .
\end{equation}
Actually, the fact that the right group and left group structures in a digroup share
a bar-unit is enough to ensure that they share all bar-units.

\begin{lemma}
\lemlabel{units}
Let $(G,\lp,\rp)$ be a dimonoid. Then $e\in G$ is a left neutral element
for $(G,\lp)$ if and only if it is a right neutral element of $(G,\rp)$.
\end{lemma}

\begin{proof}
Suppose $e\in G$ satisfies $e\lp x = x$ for all $x\in G$. Then
$1 = e\lp e\iv = e\iv$. Thus $1 = 1 \rp e$, and so 
$x \rp e = (x \rp 1) \rp e = x \rp 1 = x$. The other direction is similar.
\end{proof}

\begin{lemma}
\lemlabel{inverses}
Let $(G,\lp,\rp)$ be a digroup.
\begin{enumerate}
\item $x\lp 1 = 1\rp x$ for all $x\in G$.
\item $(x\lp y)\iv = y\iv\lp x\iv = y\iv\rp x\iv = (x\rp y)\iv$ for all $x,y\in G$.
\item $J:= \{ x\iv : x\in G \}$ is a group in which $\lp \ =\ \rp$.
\item $G\to J; x\mapsto (x\iv)\iv$ is an epimorphism of digroups
with kernel $E$.
\end{enumerate}
\end{lemma}

\begin{proof}
For (1): This follows from Lemma \lemref{left-groups}(2) and its
corresponding left group dual.

For (2): by (G3) and (G4), the inverses of $x\rp y$ and $x\lp y$ coincide,
so the result follows from Lemma \lemref{left-groups}(4) and its left
group dual.

For (3): this follows from (2), Lemma \lemref{left-groups}(5) and its
left group dual.

For (4): By Lemma \lemref{units}, the left neutral elements for $\lp$ coincide
with the right neutral elements for $\rp$, and so the desired result follows from
Lemma \lemref{left-groups}(6) and its left group dual.
\end{proof}

Incidentally, Lemma \lemref{inverses}(1) shows that one of Liu's axioms for
a digroup is redundant \cite{Liu}.

We are now ready to examine the structure of digroups in more detail.
If $(G,\lp,\rp)$ is a digroup with $E$ the set of all bar-units and $J$ the
group of all inverses, then by Lemma \lemref{left-groups}(8),
$G = J \lp E$ is isomorphic to the direct sum of $J$ and
the right zero semigroup $E$. On the other hand, by the left group
dual of Lemma \lemref{left-groups}(8), $G = E\rp J$ is isomorphic to
the direct sum of $J$ and the left zero semigroup $E$. However,
the projection of an element onto $E$ with respect to $\lp$
may not agree with the projection onto $E$ with respect to $\rp$.

\begin{example}
\exlabel{easy2}
For the digroup $G = M\times H$ of Example \exref{easy}, we have
$E = M\times \{1\}$ and $J = \{e\}\times H$. For $u\in M$, $h\in H$,
$(u,h) = (e,h) \lp (h\iv u, 1) = (u,1)\rp (e,h)$. Thus the projection of
$(u,h)$ onto $E$ with respect to $\lp$ is $(h\iv u,1)$, 
while the projection with respect to $\rp$ is $(u,1)$.
\end{example}

The appropriate generalization of conjugation to digroups is
the following. For $x$ in a digroup $G$, define
\begin{equation}
\eqlabel{digroup-conj}
x\circ y := x \lp y \rp x\iv
\end{equation}

\begin{lemma}
\lemlabel{conj-rack1}
Let $G$ be a digroup, let $\circ$ be defined by \peqref{digroup-conj},
let $E$ be the set of all bar-units, and let $J$ be the group of all inverses.
The following hold.
\begin{enumerate}
\item $x \circ (y\circ z) = (x\lp y)  \circ z = (x\rp y)\circ z$ for all $x,y,z\in G$.
\item $1\circ x = x$ and $x\circ 1 = 1$ for all $x\in G$.
\item $x\circ u\in E$ for all $x\in G$, $u\in E$.
\item $J$ acts on $E$ via $\circ$.
\end{enumerate}
\end{lemma}

\begin{proof}
For (1), we use Lemma \lemref{inverses}(2):
\begin{align*}
x\lp (y\lp z\rp y\iv)\rp x\iv &= (x\lp y)\lp z\rp (y\iv \rp x\iv)\\
&= (x\lp y)\lp z\rp (x\lp y)\iv = (x\lp y) \circ z.
\end{align*}
On the other hand, by (G3) and (G4),
\begin{align*}
(x\lp y)\circ z &= L_{\lp}(x\lp y)R_{\rp}(x\lp y)\iv z = L_{\lp}(x\rp y)R_{\rp}(x\rp y)\iv z \\
&= (x\rp y)\circ z.
\end{align*}
For (2): $1\circ x = 1\lp x\rp 1 = x$ and $x\circ 1 = (x \lp 1)\rp x\iv = (1\rp x)\rp x\iv = 1$,
by Lemma \lemref{inverses}(1).

\noindent For (3): By (G4), for $x\in G$, $u\in E$, we have
$(x\circ u)\lp y = (x\lp u\rp x\iv)\lp y = x\lp u \lp x\iv \lp y = 
x\lp x\iv \lp y = 1\lp y = y$. Similarly (or by Lemma \lemref{units}), 
$y\rp (x\circ u) = y$. Thus $x\circ u\in E$.

\noindent Finally, (4) follows from (1), (2), (3), and Lemma \lemref{inverses}(3).
\end{proof}

We now turn to this section's main result, which shows that
every digroup has the form of Example \exref{easy}.

\begin{theorem}
Let $(G,\lp,\rp)$ be a digroup with $E\subset G$ the set of bar-units
and $J\leq G$ the group of inverses. Then $G$ is isomorphic to the
digroup $(E\times J,\lp,\rp)$ where $\lp$ and $\rp$ are defined by
\begin{eqnarray}
(u,h) \lp (v,k) &= ( h\circ v, h\lp k ) \eqlabel{lpA} \\
(u,h) \rp (v,k) &= ( u, h\rp k) \eqlabel{rpA}
\end{eqnarray}
\end{theorem}

\begin{proof}
As noted above, $G = E\rp J$ is isomorphic as a
left group to $(E\times J,\rp)$ where $\rp$ is defined by \peqref{rpA}.
Denote the isomorphism by $\theta : E\times J\to E\rp J ; (u,h) \mapsto u\rp h$.
Then
\begin{align*}
\theta(u,h) \lp \theta(v,k) &= (u\rp h)\lp (v\rp k) = u\lp h\lp v \rp k \\
&= (h\circ v) \rp h \rp k = (h\circ v) \rp (h\lp k) = \theta ((u,h)\lp (v,k))
\end{align*}
using \peqref{lpA}
\end{proof}

\section{Digroups as racks and their Leibniz algebras}
\seclabel{digroup-rack}

We have already seen at the end of the last section how
conjugation in groups generalizes to digroups. Now we
draw further connections between this and racks.

\begin{lemma}
\lemlabel{conj-rack}
Let $G$ be a digroup and let $\circ$ be defined by \peqref{digroup-conj}.
The following properties hold for all $x,y,z\in G$.
\begin{enumerate}
\item $x\lp y = (x\circ y)\lp x$
\item $x \circ (y\lp z) = (x \circ y)\lp (x\circ z)$
\item $x \circ (y\rp z) = (x \circ y)\rp (x\circ z)$
\item $x \circ (y\circ z) = (x\circ y) \circ (x\circ z)$
\end{enumerate}
\end{lemma}

\begin{proof}
For (1), we use (G3):
$$
x\lp y = (x\lp y\rp x\iv) \rp x = (x\lp y\rp x\iv)\lp x = (x\circ y)\lp x.
$$
For (2), we use (G3):
\begin{align*}
x\lp (y\lp z)\rp x\iv &= x\lp y\lp x\iv\lp x\lp z \rp x\iv = (x\lp y\lp x\iv)\lp (x\circ z) \\
&= (x\lp y\rp x\iv)\lp (x\circ z) = (x\circ y)\lp (x\circ z).
\end{align*}
For (3), we use (G4): 
\begin{align*}
x\lp (y\rp z)\rp x\iv &= x\lp y\rp x\iv\rp x\rp z \rp x\iv = (x\circ y)\rp (x\rp z\rp x\iv) \\
&= (x\circ y)\rp (x\lp z\rp x\iv) = (x\circ y)\rp (x\circ z).
\end{align*}
For (4), we apply Lemma \lemref{conj-rack1}(1) twice:
$$
x\circ (y\circ z) = (x\lp y)\circ z = ((x\circ y)\lp x)\circ z = (x\circ y)\circ (x\circ z)
$$
\end{proof}

\begin{corollary}
\corolabel{digroup-rack}
Let $G$ be a digroup, and let $\circ$ be defined by \peqref{digroup-conj}.
Then $(G,\circ,1)$ is a rack.
\end{corollary}

As before, we are interested in this situation when the digroup is a manifold.

\begin{Definition}
\deflabel{lie-digroup}
A \emph{Lie digroup} $(G,\lp,\rp)$ is a smooth manifold $G$
with the structure of a digroup such that the digroup operations
$\lp,\rp : G\times G\to G$ and the inversion $(\cdot)\iv : G\to G$
are smooth mappings.
\end{Definition}

The following is an immediate consequence of the definitions
and of Theorem \thmref{rack-Leibniz}.

\begin{lemma}
\lemlabel{induced}
If $(G,\lp,\rp)$ is a Lie digroup, then its induced rack is
a Lie rack. Thus the tangent space $T_1G$ has the structure
of a Leibniz algebra. 
\end{lemma}

\begin{example}
\exlabel{lie-digroup}
Let $H$ be a Lie group, $V$ an $H$-module, and set
$G = V\times H$. As in Example \exref{easy}, define 
\begin{eqnarray}
(u,A) \lp (v,B) &= (Av, AB) \eqlabel{lp}\\
(u,A) \rp (v,B) &= (u, AB) \eqlabel{rp}
\end{eqnarray}
for all $u,v\in V$, $A,B\in H$.
Then $(G,\lp,\rp)$ is a Lie digroup, which we call a
\emph{linear Lie digroup}. The distinguished bar-unit
is $(0,1)$. The inverse of $(u,A)$ is $(0,A\iv)$.
\end{example}

Now we examine the induced rack structure of the linear Lie digroups.

\begin{theorem}
\thmlabel{main}
Let $G = V\times H$ with $H$ a Lie group and $V$ an $H$-module.
If $(G,\lp,\rp)$ is the linear Lie digroup structure defined by
\peqref{lp}-\peqref{rp}, then the induced rack $(G,\circ,\mathbf{1})$
is the linear Lie rack defined by \peqref{lie-rack}. Conversely, every
linear Lie rack is induced from a linear Lie digroup.
\end{theorem}

\begin{proof}
This all follows from a calculation:
$$
(u,A)\circ (v,B) = (u,A)\lp (v,B)\rp (0,A\iv) = (Av,ABA\iv)
$$
for all $u,v\in V$, $A,B\in H$.
\end{proof}

Finally, we obtain Lie's third theorem for split Leibniz algebras.

\begin{corollary}
\corolabel{main}
Let $\fg$ be a split Leibniz algebra. Then there exists a 
linear Lie digroup $G$ with tangent Leibniz algebra isomorphic
to $\fg$.
\end{corollary}

\begin{proof}
By Theorem \thmref{split-rack}, there exists a linear Lie rack
$G$ with tangent algebra isomorphic to $\fg$. By Theorem \thmref{main},
$G$ has a linear Lie digroup structure that induces the rack structure.
\end{proof}

\section{Conclusion}
\seclabel{conclusion}

Of course, not every Leibniz algebra splits, because not every $\fh$-module
epimorphism $\fg \to \fh \to 0$ splits. A class of nonsplit Leibniz algebras
frequently mentioned in the literature \cite{KS2} can be described as follows:
let $(\fg, \lie{\cdot}{\cdot})$ be a Lie algebra, and let $D\in \Der(\fg)$
satisfy $D^2 = 0$. Define a new bracket by
$\lie{X}{Y}_D := \lie{X}{DY}$ for all $X,Y\in \fg$. Then
$(\fg, \lie{\cdot}{\cdot}_D)$ is a Leibniz algebra. In general, $\fg$
does not split over any ideal lying between
$\cS = \langle \lie{X}{DX} : X \in \fg \rangle$ and
$\ker(\ad) = \{ Y : \lie{X}{DY} = 0\}$.

Thus the general coquecigrue problem remains open. Digroups provide
a partial solution, and indicate that the correct notion of coquecigrue for
an arbitrary Leibniz algebra should include the notion of digroup as a
special case. It is also reasonable to conjecture that whatever the 
coquecigrue turns out to be, it should induce its corresponding Leibniz
algebra via an associated Lie rack.


\references
\lastpage
\end{document}